\documentclass{article}
\begin{document}
 
\title{\bf Extended Hilbert's Nullstellensatz}
\author{Keqin Liu\\Department of Mathematics\\The University of British Columbia\\Vancouver, BC\\
Canada, V6T 1Z2}
\date{March 17, 2006}
\maketitle

\begin{abstract} We prove the extended Hilbert's Nullstellensatz in the context of \linebreak 
Hu-Liu polynomial trirings.
\end{abstract}

Which kinds of noncommutative rings are suitable for extending algebra geometry? Different attempts have been made to answer this question, but a satisfactory answer is still in hiding. My attempt at answering this question comes from the trivial extension of a ring by a bimodule over the ring. The trivial extension $R$ $\triangleright$$<$ $S$ of a ring $R$ by a $R$-bimodule $_{_R}S_{_R}$ has been used in both algebra geometry and commutative algebras for a long time (\cite{F} and \cite{M}). Even the $R$-bimodule $_{_R}S_{_R}$ is itself a ring, the multiplicative structure on the ring $_{_R}S_{_R}$ does not play any role in the trivial extension 
$R$ $\triangleright$$<$ $S$, and the researchers who make use of the trivial extension 
$R$ $\triangleright$$<$ $S$ have not paid attension to the multiplicative structure on the ring $_{_R}S_{_R}$. Simply speaking, my idea of choosing a class of noncommutative rings is not to forget the multiplicative structure on the ring $_{_R}S_{_R}$ while using the trivial extension $R$ $\triangleright$$<$ $S$. If we combine the ring product on $S$ with the bimodule actions on $_{_R}S_{_R}$ by using the Hu-Liu triassociative law, then we get a triring structure on the 
trivial extension \linebreak $R$ $\triangleright$$<$ $S$, which was introduced in \cite{Liu3}. In particular, if $R$ is a commutative ring and the $R$-bimodule $_{_R}S_{_R}$ is also a commutative ring, then the resulting triring on 
$R$ $\triangleright$$<$ $S$ is called a Hu-Liu triring. A Hu-Liu triring 
$R$ $\triangleright$$<$$S$ is a noncommutative ring with respect to the ring product on the trivial extension $R$ $\triangleright$$<$ $S$ if the left $R$-module $_{_R}S$ is different from the right $R$-module $S_{_R}$. These kinds of Hu-Liu trirings are a class of noncommutative rings which are very close to commutative rings.

\medskip
Based on the curiosity to extend algebraic geometry in the context of Hu-Liu trirings, I started on the study of affine trialgebraic sets in Chapter 5 of \cite{Liu3}. This paper is the continuation of the study. The main result is the extended Hilbert's Nullstellensatz. 

\medskip
Throughout this paper, the rings considered are assumed to have an identity, and a ring $R$ is also denoted by $(\,R, \, +, \, \cdot\,)$ to emphasize that $\cdot$ is the ring product on the ring $R$.

\medskip
We begin this paper by recalling some basic definitions and facts from \cite{Liu3}. 
A ring  $(\,R, \, +, \, \cdot\,)$ is called  a 
{\bf $\mathcal{Z}_2$-ring} if there exist two subgroups $R_0$ and $R_1$ of the additive group of $R$, called the {\bf even part} and {\bf odd part} of $R$ respectively, such that 
$R=R_0\oplus R_1$ (as Abelian groups) and 
$$
R_0R_0\subseteq R_0,\quad R_0R_1+R_1R_0\subseteq R_1, \quad R_1R_1=0.
$$
A $\mathcal{Z}_2$-ring is nothing more than a trivial extension of a ring by a bimodule over the ring. Recently, $\mathcal{Z}_2$-rings or trivial extensions of a ring by a bimodule have been used to generalize the Lie correspondence, to extend the Frobenius' Theorem about the finite dimensional division real associative algebras, and to introduce new generalizations of Jordan algebras and Malcev algebras in \cite{hl2} and \cite{hl3}.

\medskip
Let $(\,R=R_0\oplus R_1 , \, + , \, \cdot \,)$ be a $\mathcal{Z}_2$-ring.
$R$ is called a {\bf triring} if  exists a  binary operation $\sharp$ on the odd part $R_1$ of $R$ such that $(\,R_1 , \, + , \, \,\sharp\, \,)$ is a ring, and the two associative products $\cdot$ and $\,\sharp\,$ satisfy the {\bf Hu-Liu  triassociative law}:
\begin{equation}
x (\alpha \,\sharp\, \beta) =(x \alpha) \,\sharp\, \beta\quad\mbox{and}\quad
(\alpha \,\sharp\, \beta ) x =\alpha \,\sharp\, (\beta x), 
\end{equation}
where $x\in R$ and $\alpha$, $\beta\in R_1$. A triring $R=R_0\oplus R_1$ is sometimes denoted by
$$(\,R , \, + , \, \cdot , \, \,\sharp\, \,)\quad\mbox{or}\quad 
(\,R=R_0\oplus R_1 , \, + , \, \cdot ,\,  \,\sharp\, \,),$$
where the associative product $\sharp$ on the odd part $R_1$ is called the {\bf local product}, and the identity $1^\sharp$ of the ring $(\, R_1 , \, + , \, \,\sharp\,\,)$ is called the {\bf local identity} of the triring $R$.  A subset $I$ of a triring 
$(\,R , \, + , \, \cdot , \, \,\sharp\, \,)$ is called a {\bf triideal} of $R$ if $I$ is a graded ideal of the 
$\mathcal{Z}_2$-ring $R=R_0\oplus R_1$, and $I_1=I\cap R_1$ is an ideal of the ring 
$(\,R_1, \, +, \, \,\sharp\, \,)$.

\medskip
In order to clear up the relation between trivial extensions and Hu-Liu trirings, we need the following

\bigskip
\noindent
{\bf Definition 1} Let $R$ be a ring. A $R$-bimodule $_{_R}S_{_R}$ is called a {\bf sharp bimodule} over $R$ or a {\bf sharp $R$-bimodule} if $(\,S, \, +, \, \sharp \,)$ is a ring and (1) holds for $x\in R$ and 
$\alpha$, $\beta\in S$.

\bigskip
Trirings are nothing more than trivial extensions of a ring by a sharp bimodule over the ring. In fact, if $(\,R=R_0\oplus R_1 , \, + , \, \cdot ,\,  \,\sharp\, \,)$ is a triring, then 
\linebreak $R=R_0$ $\triangleright$$<$ $R_1$ is the trivial extension of the ring $R_0$ by the sharp $R_0$-bimodule $_{_{R_0}}(R_1){_{_{R_0}}}$. Conversely, if $R$ $\triangleright$\/$<$ $S$ is the trivial extension of a ring $R$ by a sharp $R$-bimodule $_{_R}S_{_R}$, then
$$(\,\mbox{$R$ $\triangleright$$<$ $S$}=(\mbox{$R$ $\triangleright$$<$ $S$})_0\oplus 
(\mbox{$R$ $\triangleright$$<$ $S$})_1 , \, + , \, \cdot ,\,  \,\sharp\, \,)$$
is a triring, where the local product $\sharp$ is the ring product carried by the ring 
$(\,S, \, +, \, \sharp \,)$, and the even part $(\mbox{$R$ $\triangleright$$<$ $S$})_0$ and the odd part $(\mbox{$R$ $\triangleright$$<$ $S$})_1$ are given by
$$(\mbox{$R$ $\triangleright$$<$ $S$})_0:=(\, R, \, 0\,) \quad{and}\quad
(\mbox{$R$ $\triangleright$$<$ $S$})_1:=(\, 0, \, S\,).$$

\medskip
We now recall the following definition from  \cite{Liu3}.

\bigskip
\noindent
{\bf Definition 2} Let $(\,R=R_0\oplus R_1 , \, + , \, \cdot , \, \,\sharp\, \,)$ be a triring.
\begin{description}
\item[(i)] $R$ is called a {\bf commutative triring} if both $(\,R , \, + , \, \cdot \,)$ and $(\,R_1 , \, + ,\, \,\sharp\, \,)$ are commutative rings.
\item[(ii)] $R$ is called a {\bf Hu-Liu triring} if both $(\,R_0 , \, + , \, \cdot \,)$ and $(\,R_1 , \, + ,\, \,\sharp\, \,)$ are commutative rings.
\item[(iii)] $R$ is called an {\bf algebraic closed $3$-trifield} if both 
$(\,R_0, \, +, \, \cdot \,)$ and 
\linebreak $(\,R_1, \, +, \, \,\sharp\, \,)$ are algebraic closed fields, and 
$R_1=R_0 1^\sharp =1^\sharp R_0$, where $1^\sharp$ is the local identity of the triring $R$.
\end{description}

\bigskip
Commutative trirings are a class of commutative rings with many zero divisors. Since the zero divisors in a commutative triring are controlled nicely by the local product, the theory of commutative algebra can be extended satisfactorily in the context of commutative trirings. For example, the book \cite{S} can be rewritten in the context of commutative trirings.

\medskip
Let $(\,\mathbf{k}=\mathbf{k}_0\oplus \mathbf{k}_1, \, +,\, \cdot, \, \,\sharp\, \,)$ be an algebraic closed $3$-trifield with a local identity $1^{\sharp}$, where $\mathbf{k}_0$ and $\mathbf{k}_1$ are the even part and the odd part of $\mathbf{k}$ respectively. Let $\mathbf{k}^{\sharp}[n]:=\mathbf{k}^{\sharp} [x_{(1)}, x_{(2)}, \cdots , x_{(n)}]$ be the Hu-Liu polynomial triring\footnote{See Section 4.1 in \cite{Liu3} for the construction of the Hu-Liu polynomial triring.} over $\mathbf{k}$ in $n$ indeterminates $x_{(1)}$, $x_{(2)}$, $\cdots$, $x_{(n)}$. The Hu-Liu polynomial triring 
$\mathbf{k}^{\sharp}[n]=\mathbf{k}^{\sharp} [x_{(1)}, x_{(2)}, \cdots , x_{(n)}]$ is a Hu-Liu triring, its even part 
$$\mathbf{k}^{\sharp}[n]_0=\mathbf{k}_0 [x_{(1)0}, x_{(2)0}, \cdots , x_{(n)0}]$$
is the polynomial ring over the algebraic closed field $(\,\mathbf{k}_0, \, +,\, \cdot \,)$ in $n$ indeterminates $x_{(1)0}$, $x_{(2)0}$, $\cdots$, $x_{(n)0}$, and its odd part
$$
\mathbf{k}^{\sharp}[n]_1=\mathbf{k}_1 [x_{(1)0}1^\sharp, \cdots , x_{(n)0}1^\sharp, x_{(1)1}, 
\cdots , x_{(n)1}, 1^\sharp x_{(1)0}, \cdots , 1^\sharp x_{(n)0}]
$$
is the polynomial ring over the algebraic closed field $(\,\mathbf{k}_1, \, +,\, \,\sharp\, \,)$ in $3n$ indeterminates $x_{(1)0}1^\sharp$, $\cdots$, $x_{(n)0}1^\sharp$, $x_{(1)1}$, $\cdots$ , $x_{(n)1}$, $1^\sharp x_{(1)0}$, $\cdots$, $1^\sharp x_{(n)0}$. 
Let $\mathbf{A}^n_\mathbf{k}$ be the affine $n$-trispace over $\mathbf{k}$, where $\mathbf{A}^n_\mathbf{k}$ is defined by
$$
\mathbf{A}^n_\mathbf{k}:=\{\, (a_{(1)}, \,a_{(2)} , \, \cdots , \,a_{(n)}) \,| \, \mbox{$a_{(j)}\in \mathbf{k}$ for $1\le j\le n$} \,\}.
$$ 
If $a=(a_{(1)}, \,a_{(2)} , \, \cdots , \,a_{(n)})\in \mathbf{A}^n_\mathbf{k}$, we  define
\begin{eqnarray*}
a_0:&=&(a_{(1)0}, \,a_{(2)0} , \, \cdots , \,a_{(n)0} )\in \mathbf{A}^n_{\mathbf{k}_0}, \\
a_01^{\sharp}:&=&(a_{(1)0}1^{\sharp}, \,a_{(2)0}1^{\sharp} , \, \cdots , \,a_{(n)0}1^{\sharp}),\\ 
a_1:&=&(a_{(1)1}, \,a_{(2)1} , \, \cdots , \,a_{(n)1}),\\
1^{\sharp}a_0:&=&(1^{\sharp}a_{(1)0}, \,1^{\sharp}a_{(2)0} , \, \cdots , \,1^{\sharp}a_{(n)0} ), 
\end{eqnarray*}
where $a_{(j)}=a_{(j)0}+a_{(j)1}$, $a_{(j)0}\in\mathbf{k}_0$ and $a_{(j)1}\in\mathbf{k}_1$ are the even component and the odd component of $a_{(j)}$ for $1\le j\le n$.
Let $J=J_0\oplus J_1$ be a triideal of $\mathbf{k}^\sharp[n]$. The even affine trialgebraic set $V^\sharp _0(J)$ and the odd affine trialgebraic set $V^\sharp _1(J)$ determined by the triideal $J$ are defined by
$$
V^\sharp _0(J):=\{\, a\in \mathbf{A}^n_\mathbf{k} \, | \, 
\mbox{$F_0(a_0)=0$ for all $F_0\in J_0$} \, \}
$$
and
$$
V^\sharp _1(J):=\{\, a\in \mathbf{A}^n_\mathbf{k} \, | \, 
\mbox{$F_1(a_01^\sharp, \, a_1, \, 1^\sharp a_0)=0$ for all $F_1\in J_1$} \, \}.
$$
Clearly, $V^\sharp _1(J)\subseteq V^\sharp _0(J)$. Other basic properties of even and odd affine trialgebraic sets are given in Proposition 5.1 of \cite{Liu3}.

\medskip
Let $V_0$ and $V_1$ be a pair of subsets of the affine $n$-trispace $\mathbf{A}^n_\mathbf{k}$, where $V_1\subseteq V_0$. Define $I^\sharp (V_0, V_1)$ by
$$
I^\sharp (V_0, V_1):=\left\{\, F_0+F_1 \,\left |\, \begin{array}{c}\mbox{$F_0\in\mathbf{k}^{\sharp}[n]_0$, $F_1\in\mathbf{k}^{\sharp}[n]_1$,}\\ \mbox{$F_0(a_0)=0$ for all $a\in V_0$, and}\\
\mbox{$F_1(b_01^\sharp, \, b_1, \, 1^\sharp b_0)=0$ for all $b\in V_1$}\end{array}\right.\,\right\}.
$$
Then $I^\sharp (V_0, V_1)$ is a triideal of $\mathbf{k}^{\sharp}[n]$.

\medskip
The next proposition gives the main result of this paper.

\bigskip
\noindent
{\bf Proposition} ({\bf The Extended Hilbert's Nullstellensatz}) Let 
$\mathbf{k}=\mathbf{k}_0\oplus \mathbf{k}_1$ be an algebraic closed $3$-trifield. If $J$ is a triideal of the Hu-Liu polynomial triring 
$\mathbf{k}^{\sharp}[n]=\mathbf{k}^{\sharp} [x_{(1)}, x_{(2)}, \cdots , x_{(n)}]$, then
$I^\sharp (V^\sharp _0(J), V^\sharp _1(J))=\sqrt[\sharp]{J}$, where 
$$
\sqrt[\sharp]{J}:=\left\{\, G_0+G_1 \,\left |\, \begin{array}{c}\mbox{$G_0\in\mathbf{k}^{\sharp}[n]_0$, $G_1\in\mathbf{k}^{\sharp}[n]_1$, and}\\
\mbox{$G_0^m\in J_0$ and $G_1^{\,\sharp\, n}\in J_1$ for some $m$, $n\in\mathcal{Z}_{>0}$}\end{array}\right.\,\right\}
$$
is the graded nilradical of $J$, and $G_1^{\sharp n} :=
\underbrace{G_1 \,\sharp\, G_1  \,\sharp\, \cdots  \,\sharp\, G_1}_n$. 

\medskip
\noindent
{\bf Proof} Let $G_0+G_1$ be an element of $\sqrt[\sharp]{J}=\sqrt{J_0}\oplus \sqrt{J_1}$, where 
$G_0\in \sqrt{J_0}$ and $G_1\in \sqrt{J_1}$. Then $G_0^m\in J_0$ and $G_1^{\,\sharp\, n}\in J_1$ for some $m$, $n\in\mathcal{Z}_{>0}$. Hence, we have
$$
\Big(G_0(a_0)\Big)^m=G_0^m(a_0)=0 \quad\mbox{for all $a\in V^\sharp _0(J)$}
$$
and
$$\Big(G_1(b_01^\sharp, \, b_1, \, 1^\sharp b_0)\Big)^{\,\sharp\, n}=
G_1^{\,\sharp\, n}(b_01^\sharp, \, b_1, \, 1^\sharp b_0)=0 \quad\mbox{for all $b\in V^\sharp _1(J)$}.
$$
It follows that $G_0(a_0)=0$ for all $a\in V^\sharp _0(J)$ and 
$G_1(b_01^\sharp, \, b_1, \, 1^\sharp b_0)=0$ for all $b\in V^\sharp _1(J)$. Thus, 
$G_0+G_1\in I^\sharp (V^\sharp _0(J), V^\sharp _1(J))$. This proves 
\begin{equation}\label{eq2}
\sqrt[\sharp]{J}\subseteq I^\sharp (V^\sharp _0(J), V^\sharp _1(J)).
\end{equation}

\medskip
Conversely, let $F=F_0+F_1\in I^\sharp (V^\sharp _0(J), V^\sharp _1(J))$ with $F_i\in 
\mathbf{k}^{\sharp}[n]_i$ for $i=0$ and $1$, then 
\begin{equation}\label{eq3}
F_0(a_0)=0 \quad\mbox{for all $a_0\in V(J_0)\subseteq V^\sharp _0(J)$}
\end{equation}
and 
\begin{equation}\label{eq4}
F_1(b_01^\sharp, \, b_1, \, 1^\sharp b_0)=0 \quad\mbox{for all $b\in V^\sharp _1(J)$},
\end{equation}
where $V(J_0)$ is the ordinary affine algebraic set in the affine 
$n$-space $\mathbf{A}^n_{\mathbf{k}_0}$. By (\ref{eq3}), we have $F_0\in I(V(J_0))$, where 
$I(V(J_0))$ is the ideal determined by $V(J_0)$. Using Hilbert's Nullstellensats in the 
polynomial ring 
$$\mathbf{k}^{\sharp}[n]_0=\mathbf{k}_0 [x_{(1)0}, x_{(2)0}, \cdots , x_{(n)0}],$$ 
we get
\begin{equation}\label{eq5}
F_0\in I(V(J_0))=\sqrt {J_0}\subseteq \sqrt[\sharp]{J}.
\end{equation}

If $F_1=0$, then $F=F_0\in \sqrt[\sharp]{J}$ by (\ref{eq5}). Hence, we get
\begin{equation}\label{eq6}
\mbox{$F=F_0+F_1\in I^\sharp (V^\sharp _0(J), V^\sharp _1(J))$ and $F_1= 0$}\Longrightarrow F\in \sqrt[\sharp]{J}.
\end{equation}

If $F_1\ne 0$, we consider the Hu-Liu polynomial triring 
$$
R:=\mathbf{k}^{\sharp}[n+1]=\Big(\mathbf{k}^{\sharp}[n]\Big)^{\sharp}[x_{(n+1)}]=
\mathbf{k}^{\sharp} [x_{(1)}, \cdots , x_{(n)}, x_{(n+1)}]
$$ 
over $\mathbf{k}$ in $n+1$ indeterminates $x_{(1)}$, $\cdots$, $x_{(n)}$, $x_{(n+1)}$. Let
$$U:=J_0R_0\oplus \Big(J_1\,\sharp\, R_1+(x_{(n+1)1}\,\sharp\, F_1 -1^\sharp)\,\sharp\, R_1\Big),$$
where $R_0$ and $R_1$ are the even part and the odd parts of $R$ respectively. Clearly, $U$ is a triideal of $R$ and $J=J_0\oplus J_1\subseteq U$. If 
$$U_1=J_1\,\sharp\, R_1+(x_{(n+1)1}\,\sharp\, F_1 -1^\sharp)\,\sharp\, R_1\ne R_1,$$ 
then $V^\sharp _1(U)\ne \emptyset$ by Proposition 5.2 (ii) in \cite{Liu3}. Thus, there would exist 
$$c=(c_{(1)}, \cdots , c_{(n)}, c_{(n+1)})\in \mathbf{A}^{n+1}_\mathbf{k}$$ such that
\begin{equation}\label{eq7}
\mathbf{g}_1(c_01^{\sharp}, \, c_1, \,1^{\sharp}c_0)=0\quad\mbox{for all $\mathbf{g}_1\in U_1$,}
\end{equation}
which would imply that
\begin{equation}\label{eq8}
\mathbf{f}_1(d_01^{\sharp}, \, d_1,\,1^{\sharp}d_0)=0\quad\mbox{for all $\mathbf{f}_1\in J_1$,}
\end{equation}
where $d=(c_{(1)}, \cdots , c_{(n)})\in \mathbf{A}^{n}_\mathbf{k}$. It follows from (\ref{eq8}) that
\begin{equation}\label{eq9}
d=(c_{(1)}, \cdots , c_{(n)})\in V^\sharp _1(J).
\end{equation}
By (\ref{eq4}), (\ref{eq7}) and (\ref{eq9}), we get
\begin{eqnarray*}
0&=&(x_{(n+1)1}\,\sharp\, F_1 -1^\sharp)(c)=
c_{(n+1)1}\,\sharp\, F_1(c_01^{\sharp}, \, c_1, \,1^{\sharp}c_0) -1^\sharp\\
&=&c_{(n+1)1}\,\sharp\, F_1(d_01^{\sharp}, \, d_1, \,1^{\sharp}d_0) -1^\sharp
=c_{(n+1)1}\,\sharp\, 0-1^\sharp=-1^\sharp ,
\end{eqnarray*}
which is impossible. This proves that 
$$U_1=J_1\,\sharp\, R_1+(x_{(n+1)1}\,\sharp\, F_1 -1^\sharp)\,\sharp\, R_1= R_1 .$$
Therefore, the local identity $1^\sharp$ can be expressed as 
\begin{equation}\label{eq10}
1^\sharp=\displaystyle\sum_{j=1}^m F_{(j)}\,\sharp\, \mathbf{g}_{(j)}+
(x_{(n+1)1}\,\sharp\, F_1 -1^\sharp)\,\sharp\, \mathbf{h},
\end{equation}
where $m\in\mathcal{Z}_{>0}$, $F_{(j)}\in J_1$, $\mathbf{g}_{(j)}\in R_1$, $\mathbf{h}\in R_1$ and $1\le j\le m$. Let
$$
\mathbf{k}_1(3n):=\mathbf{k}_1\Big(x_{(1)0}1^\sharp, \cdots , x_{(n)0}1^\sharp, x_{(1)1}, 
\cdots , x_{(n)1}, 1^\sharp x_{(1)0}, \cdots , 1^\sharp x_{(n)0}\Big)
$$
be the fractional field of the polynomial ring $\mathbf{k}_1[3n]$ over the field 
$(\,\mathbf{k}_1, \, +,\, \,\sharp\, \,)$ in $3n$ indeterminates $x_{(1)0}1^\sharp$, $\cdots$, $x_{(n)0}1^\sharp$, $x_{(1)1}$, $\cdots$ , $x_{(n)1}$, $1^\sharp x_{(1)0}$, $\cdots$, $1^\sharp x_{(n)0}$, where
$$
\mathbf{k}_1[3n]:=\mathbf{k}_1\Big[x_{(1)0}1^\sharp, \cdots , x_{(n)0}1^\sharp, x_{(1)1}, 
\cdots , x_{(n)1}, 1^\sharp x_{(1)0}, \cdots , 1^\sharp x_{(n)0}\Big].
$$
By the universal property of the polynomial ring 
$$
\mathbf{k}_1[3(n+1)]:=\Big(\mathbf{k}_1[3n]\Big)\Big[x_{(n+1)0}1^\sharp , x_{(n+1)1}, 
1^\sharp x_{(n+1)0}\Big],
$$
there exists a ring homomorphism $\psi : \mathbf{k}_1[3(n+1)]\to \mathbf{k}_1(3n)$ such that
$$
\psi \left|\Big(\mathbf{k}_1[3n]\Big)\right.=\mbox{ the identity map}
$$
and
$$
\psi\Big(x_{(n+1)0}1^\sharp\Big)=0=\psi\Big(1^\sharp x_{(n+1)0}\Big), \quad 
\psi\Big(x_{(n+1)1}\Big)=\frac{1}{F_1}.
$$

After applying $\psi$ to (\ref{eq10}), we get
\begin{eqnarray}
&&1^\sharp=\psi(1^\sharp)
=\displaystyle\sum_{j=1}^m \psi\left(F_{(j)}\right)\,\sharp\, \psi\left(\mathbf{g}_{(j)}\Big(x_{(1)0}1^\sharp, \cdots , x_{(n)0}1^\sharp, x_{(n+1)0}1^\sharp, \right.\nonumber\\
&&\quad \left. x_{(1)1}, \cdots , x_{(n)1}, x_{(n+1)1}, 1^\sharp x_{(1)0}, \cdots , 1^\sharp x_{(n)0}, 1^\sharp x_{(n+1)0}\Big)\right)\nonumber\\
&=&\displaystyle\sum_{j=1}^m F_{(j)}\,\sharp\, \mathbf{g}_{(j)}\Big(x_{(1)0}1^\sharp, \cdots , x_{(n)0}1^\sharp, 0, x_{(1)1}, \cdots , x_{(n)1}, \frac{1}{F_1},\nonumber\\
\label{eq11}
&&\quad  1^\sharp x_{(1)0}, \cdots , 1^\sharp x_{(n)0}, 0\Big).
\end{eqnarray}
Let $s_j$ be the degree of 
$$
\mathbf{g}_{(j)}\Big(x_{(1)0}1^\sharp, \cdots , x_{(n)0}1^\sharp, 0, x_{(1)1}, \cdots , x_{(n)1}, \frac{1}{F_1},  1^\sharp x_{(1)0}, \cdots , 1^\sharp x_{(n)0}, 0\Big)
$$
when considered as a polynomial in $\displaystyle\frac{1}{F_1}$ with coefficients in $\mathbf{k}_1[3n]$. Then
$$
(F_1)^{\sharp s}\,\sharp\,\mathbf{g}_{(j)}\Big(x_{(1)0}1^\sharp, \cdots , x_{(n)0}1^\sharp, 0, x_{(1)1}, \cdots , x_{(n)1}, \frac{1}{F_1},  1^\sharp x_{(1)0}, \cdots , 1^\sharp x_{(n)0}, 0\Big)
$$
is in $\mathbf{k}_1[3n]$, where $s:=\max\{\, s_j \,|\, 1\le j\le m\,\}$. Using this fact and multiplying both sides of (\ref{eq11}) by $(F_1)^{\sharp s}$, we get 
\begin{equation}\label{eq12}
(F_1)^{\sharp s}\in \displaystyle\sum_{j=1}^m F_{(j)}\,\sharp\, \left(\mathbf{k}_1[3n]\right)\subseteq J_1.
\end{equation}
By (\ref{eq5}) and (\ref{eq12}), we have
\begin{equation}\label{eq13}
\mbox{$F=F_0+F_1\in I^\sharp (V^\sharp _0(J), V^\sharp _1(J))$ and $F_1\ne 0$}\Longrightarrow F\in \sqrt[\sharp]{J}.
\end{equation}
It follows from (\ref{eq6}) and (\ref{eq13}) that
\begin{equation}\label{eq14}
\sqrt[\sharp]{J}\supseteq I^\sharp (V^\sharp _0(J), V^\sharp _1(J)).
\end{equation}
The extended Hilbert's Nullstellensatz now follows from (\ref{eq2}) and (\ref{eq14}).

\hfill\raisebox{1mm}{\framebox[2mm]{}}


\bigskip
\begin{thebibliography}{99}
\bibitem{F} Robert M. Fossum, Phillip A. Griffith and Idun Reiten, \textsl{Trivial Extensions of Abelian Categories}, \quad Lecture Notes in Mathematics, Vol.456, Springer-Verlag, Berlin, 1975.
\bibitem{Liu3} Keqin Liu, \textsl{Introduction to Trirings}, \quad Research Monographs in Mathematics {\bf 2}, 153 Publishing, 2006.
\bibitem{hl2} Keqin Liu, \textsl{Number-like objects and the extended Lie correspondence}, arXiv: math.RA/0604179 v1 8 Apr 2006.
\bibitem{hl3} Keqin Liu, \textsl{Generalizations of Jordan algebras and Malcev algebras}, arXiv: math.RA/0606628 v1 24 Jun 2006.
\bibitem{M} Hideyuki Matsumura, \textsl{Commutative Algebra ({\small{second edition}})}, The Benjamin/Cummings Publishing Company, Inc., 1980.
\bibitem{S} Rodney Y. Sharp, \textsl{Steps in Commutative Algebra}, London Mathematical Society Student Texts {\bf 19}, Cambridge University Press, 1990.
\end{thebibliography}
\end{document}